\documentclass{amsart}
\usepackage{amssymb}
\usepackage{graphicx}
\usepackage{amsmath}

\newtheorem{theorem}{Theorem}

\newtheorem{corollary}[theorem]{Corollary}

\newtheorem{lemma}[theorem]{Lemma}

\newtheorem{proposition}[theorem]{Proposition}

\newcommand{\lm}{\operatorname{lim}}
\begin{document}

\title{Products of compact filters and applications to classical product
theorems}
\author{Fr\'{e}d\'{e}ric Mynard}
\address{Department of Mathematical Sciences, Georgia Southern University, Statesboro, GA 30460-8093}
\email{fmynard@georgiasouthern.edu}
\maketitle

\begin{abstract}
Two results on product of compact filters are shown to be the common
principle behind a surprisingly large number of theorems.
\end{abstract}

\section{Introduction}

The terminology and notations are those of the companion paper \cite%
{myn.relations}. In particular, two families $\mathcal{A}$ and $\mathcal{B}$
of subsets of $X$ are said to \emph{mesh,} in symbol $\mathcal{A}\#\mathcal{B%
},$ if $A\cap B\neq \varnothing $ whenever $A\in \mathcal{A}$ and $B\in 
\mathcal{B}.$ Given a class $\mathbb{D}$ of filters on $X$ and $A\subset X,$
we call a filter $\mathcal{F}$ (on $X$) $\mathbb{D}$-\emph{compact at} $A$
if 
\begin{equation*}
\mathcal{D}\in \mathbb{D},\mathcal{D}\#\mathcal{F}\Longrightarrow \mathrm{%
adh}\mathcal{D}\cap A\neq \varnothing .
\end{equation*}%
The context of the present paper is that of \emph{convergence spaces }as
defined in \cite{myn.relations} and therefore $\mathrm{adh}\mathcal{D}$
denotes the union of limit sets of filters finer than $\mathcal{D}.$ The
notion derives from \emph{total nets} introduced by Pettis \cite{pettis} and
turned out to be very useful in a variety of contexts, for instance in \cite%
{JLM}, \cite{labuda.boundary}, \cite{LabudaDavis.inherent}, in \cite{DGL}, 
\cite{active}, \cite{D.comp}, \cite{Labuda02}, \cite{cascales} under the
name of compactoid filter, \cite{pettis}, \cite{vaughan78}, \cite{total}
under the name of total filter.

In \cite{myn.relations}, many classes of maps are characterized as relations
preserving $\mathbb{D}$-compactness of filters. The aim of this paper is to
establish a pair of theorems on product of $\mathbb{D}$-compact filters and
show that, in view of the results of \cite{myn.relations}, they are the
common principle behind a surprising number of results of stability under
product of global properties (variants of compactness), local properties (Fr%
\'{e}chetness and variants) and maps (variants of quotient and perfect maps).

\section{Characterization of $\mathbb{D}$-compact filters in terms of
products \label{sec:Dcompprod}}

The goal of this section is to show that the classical Mr\'{o}wka-Kuratowski
Theorem characterizing compactness of $X$ in terms of closedness of the
projections $p_{Y}:X\times Y\rightarrow Y$ for every topological space $Y$
and its variants for other type of compactness (e.g., countable compactness,
Lindel\"{o}fness), as well as product characterizations of various types of
maps are all instances of a simple result on $\mathbb{D}$-compact filters.

If $\mathbb{D}$ and $\mathbb{J}$ are two classes of filters, we say that $%
\mathbb{J}$ is $\mathbb{D}$-\emph{composable} if for every $X$ and $Y,$ the
(possibly degenerate) filter $\mathcal{H}\left( \mathcal{F}\right) $=$\{HF:$ 
$H\in \mathcal{H},F\in \mathcal{F}\}^{\uparrow }$ (\footnote{$HF=\{y\in
Y:(x,y)\in H$ and $x\in F\}.$}) belongs to $\mathbb{J}(Y)$ whenever $%
\mathcal{F\in }\mathbb{J}(X)$ and $\mathcal{H}\in \mathbb{D}(X\times Y),$
with the convention that every class of filters contains the degenerate
filter. If a class $\mathbb{D}$ is $\mathbb{D}$-composable, we simply say
that $\mathbb{D}$ is \emph{composable}. Notice that 
\begin{equation}
\mathcal{H}\#\left( \mathcal{F\times G}\right) \Longleftrightarrow \mathcal{H%
}\left( \mathcal{F}\right) \#\mathcal{G}\Longleftrightarrow \mathcal{H}%
^{-}\left( \mathcal{G}\right) \#\mathcal{F},  \label{eq:grill}
\end{equation}%
where $\mathcal{H}^{-}\left( \mathcal{G}\right) $=$\{H^{-}G=\{x\in
X:(x,y)\in H$ and $y\in G\}:H\in \mathcal{H},G\in \mathcal{G}\}^{\uparrow }.$

\begin{theorem}
\label{th:Dcompactproduct}Let ($X,\xi )$ be a convergence space, $A\subset
X, $ and let $\mathcal{F}$ be a filter on $X$. Let $\mathbb{D}$ be a
composable class of filters. The following are equivalent:

\begin{enumerate}
\item $\mathcal{F}$ is $\mathbb{D}$-compact at $A;$

\item For every convergence space $Y$ and every compact $\mathbb{D}$-filter $%
\mathcal{G}$ at $B\subset Y,$ the filter $\mathcal{F}\times \mathcal{G}$ is $%
\mathbb{D}$-compact at $A\times B;$

\item For every $\mathbb{D}$-based atomic (\footnote{%
A topological space with at most one non-isolated point is called \emph{%
atomic. }Such spaces have been also called \emph{point-spaces} and \emph{%
prime topological spaces}.}) topological space $Y,$ every $y\in Y$ and every 
$\mathcal{G}$ such that $y\in \lm_{Y}\mathcal{G},$ the filter $\mathcal{F}%
\times \mathcal{G}$ is $\mathbb{F}_{1}$-compact at $A\times \{y\}.$
\end{enumerate}
\end{theorem}

\begin{proof}
$\left( 1\Longrightarrow 2\right) $.

Let $\mathcal{D}\in \mathbb{D}(X\times Y)$ such that $\mathcal{D}\#\left( 
\mathcal{F}\times \mathcal{G}\right) .$ The filter $\mathcal{D}^{-}\mathcal{%
(G)}\in \mathbb{D}(X)$ because $\mathcal{G\in }\mathbb{D}(Y)$ and $\mathbb{D}
$ is composable. Moreover $\mathcal{D}^{-}\mathcal{(G)}\#\mathcal{F}$ so
that $\mathrm{adh}_{X}\mathcal{D}^{-}\mathcal{(G)}\cap A\neq \emptyset .$
Consequently, there exists a filter $\mathcal{W}$ with $x\in \lm_{X}%
\mathcal{W}\cap A$ such that $\mathcal{W}\#\mathcal{D}^{-}\mathcal{(G)}.$
Therefore $\mathcal{D(W)}\#\mathcal{G}$ and $\mathrm{adh}_{Y}\mathcal{D(W)}%
\cap B\neq \emptyset $ by compactness of $\mathcal{G}.$ In other words,
there is a filter $\mathcal{U}$ with $y\in \lm_{Y}\mathcal{U}\cap B$ such
that $\mathcal{U}\#\mathcal{D(W)}.$ Consequently, $(x,y)\in \mathrm{adh}%
_{X\times Y}\mathcal{D}$ because $\left( \mathcal{W}\times \mathcal{U}%
\right) \#\mathcal{D}.$

$\left( 2\Longrightarrow 3\right) $ is obvious.

$\left( 3\Longrightarrow 1\right) $.

Assume that $\mathcal{F}$ is not $\mathbb{D}$-compact at $A.$ Then, there
exists a $\mathbb{D}$-filter $\mathcal{D}\#\mathcal{F}$ such that $\mathrm{%
adh}_{\xi }\mathcal{D}\cap A=\emptyset .$ Chose any point $x_{0}$ in $X$ and
let $Y$ be a copy of $X$ endowed with the atomic topology $\tau $ defined by 
$\mathcal{N}_{\tau }(x_{0})=\mathcal{D}\wedge (x_{0}).$ Then $\mathcal{F}%
\times \mathcal{N}_{\tau }(x_{0})$ is \emph{not }$\mathbb{F}_{1}$-compact at 
$A\times \{x_{0}\}$. Indeed, $\{(x,x):x\neq x_{0}\}\#\mathcal{F}\times 
\mathcal{N}_{\tau }(x_{0})$ because $\mathcal{D}\#\mathcal{F},$ but $%
\mathrm{adh}_{\xi \times \tau }\{(x,x):x\neq x_{0}\}\cap \left(
A\times \{x_{0}\}\right) =\emptyset .$ For a filter on $\{(x,x):x\neq
x_{0}\} $ is of the form $\mathcal{G}\times \mathcal{G}\ $and if $x_{0}\in
\lm_{\tau }\mathcal{G}$ then $\mathcal{G}\geq \mathcal{D}$, so that $%
\lm_{\xi }\mathcal{G}\cap A=\emptyset .$
\end{proof}

A relation $R:X\rightrightarrows Y$ is called $\mathbb{D}$-\emph{compact} if 
$R\left( \mathcal{F}\right) $ is $\mathbb{D}$-compact at $RA$ whenever $%
\mathcal{F}$ is $\mathbb{D}$-compact at $A.$ As observed in \cite[section 10]%
{D.comp}, preservation of closed sets by a map $f:(X,\xi )\rightarrow
(Y,\tau )$ is equivalent to $\mathbb{F}_{1}$-compactness of the inverse map $%
f^{-}$ when $(X,\xi )$ is topological, but not if $\xi $ is a general
convergence. More precisely, calling a map $f:(X,\xi )\rightarrow (Y,\tau )$ 
\emph{adherent }\cite{D.comp} if 
\begin{equation*}
y\in \mathrm{adh}_{\tau }f(H)\Longrightarrow \mathrm{adh}%
_{\xi }H\cap f^{-}y\neq \emptyset ,
\end{equation*}%
we have:

\begin{lemma}
\begin{enumerate}
\item \label{lem:closed} A map $f:(X,\xi )\rightarrow (Y,\tau )$ is adherent
if and only if $f^{-}:(Y,\tau )\rightrightarrows (X,\xi )$ is $\mathbb{F}%
_{1} $-compact;

\item If $f:(X,\xi )\rightarrow (Y,\tau )$ is adherent, then it is closed;

\item If $f:(X,\xi )\rightarrow (Y,\tau )$ is closed and if adherence of
sets are closed in $\xi $ (in particular if $\xi $ is a topology), then $f$
is adherent.
\end{enumerate}
\end{lemma}

\cite[Theorem 13]{myn.relations} shows that a map $f:X\rightarrow Y$ is $%
\mathbb{D}$-\emph{perfect} (that is, adherent with $\mathbb{D}$-compact
fibers) if and only if the inverse map $f^{-}:Y\rightrightarrows X$ is $%
\mathbb{D}$-compact. Hence, applied for $\mathcal{F}=\{X\}=\{A\},$ Theorem %
\ref{th:Dcompactproduct} rephrases as:

\begin{corollary}
\label{cor:Dcompactproj}Let $\mathbb{D}$ be a composable class of filters
and let $X$ be a convergence space. The following are equivalent:

\begin{enumerate}
\item $X$ is $\mathbb{D}$-compact;

\item for every $\mathbb{D}$-based convergence space $Y$, the projection $%
p_{Y}:X\times Y\rightarrow Y$ is $\mathbb{D}$-perfect;

\item for every $\mathbb{D}$-based atomic topological space $Y,$ the
projection $p_{Y}:X\times Y\rightarrow Y$ is adherent.
\end{enumerate}
\end{corollary}

\begin{proof}
$(1\Longrightarrow 2)$ because the fact that $\{X\}\times \mathcal{G}$ is $%
\mathbb{D}$-compact at $X\times \{y\}$ for every $\mathbb{D}$-filter $%
\mathcal{G}$ such that $y\in \lm_{Y}\mathcal{G}$ amounts to $\mathbb{D}$%
-compactness of $p_{Y}^{-}:Y\rightrightarrows X\times Y,$ which implies $%
\mathbb{D}$-perfectness of $p_{Y}:X\times Y\rightarrow Y.$

$(2\Longrightarrow 3)$ by definition, and $(3\Longrightarrow 1)$ because if $%
p_{Y}:X\times Y\rightarrow Y$ is adherent for every $\mathbb{D}$-based
atomic topological space $Y,$ then for every topological space $Y,$ every $%
y\in Y$ and every $\mathbb{D}$-filter $\mathcal{G}$ that converges to $y,$
the filter $\{X\}\times \mathcal{G}$ is $\mathbb{F}_{1}$-compact at $X\times
\{y\}.$ In view of Theorem \ref{th:Dcompactproduct}, $\{X\}$ is compact ,
that is, $X$ is compact.
\end{proof}

In particular, for a topological space $X,$ $\mathbb{D}$-compactness amounts
to $\left( \int \mathbb{D}\right) $-compactness (\footnote{%
If $\mathcal{F}$ is a filter on $X$ and $\mathcal{G}:X\rightarrow \mathbb{F}%
X $ then the \emph{contour of} $\mathcal{G}$ \emph{along }$\mathcal{F}$ is
the filter on $X$ defined by%
\begin{equation*}
\int_{\mathcal{F}}\mathcal{G}=\bigvee_{F\in \mathcal{F}%
}\bigwedge_{x\in F}\mathcal{G}(x).
\end{equation*}%
\par
A well-capped tree with least element is called a \textit{filter cascade} if
its every (non maximal) element is a filter on the set of its immediate
successors.
\par
A map $\Phi :V\setminus \{\varnothing _{V}\}\rightarrow X$, where $V$ is a
filter cascade with the least element $\varnothing _{V}$, is called a 
\textit{multifilter on} $X$. If${\ }\mathbb{D}$ is a class of filters, we
call ${\ }\mathbb{D}$-\emph{multifilter} a multifilter with a cascade of ${\ 
}\mathbb{D}$-filters as domain. For each $v\in V,$ the subset $V(v)$ of $V$
formed by $v$ and its successors is also a cascade. The \textit{contour} of $%
\Phi :V\setminus \{\varnothing _{V}\}\rightarrow X$ is defined by induction
to the effect that $\int \Phi \ $is the filter generated by $\varnothing
_{V} $ on $\Phi (\max V)$ if $V=\{\varnothing _{V}\}\cup $ $\max V$, and%
\begin{equation*}
\int \Phi =\int_{\varnothing _{V}}\left( \int \Phi |_{V(.)}\right)
\end{equation*}%
otherwise. With each class $\mathbb{D}$ of filters we associate the class $%
\int \mathbb{D}$ of all $\mathbb{D}$-contour filters, i.~e., the contours of 
$\mathbb{D}$-multifilter. See \cite{cascades} for details.}) so that, in
view of \cite[Lemma 6]{myn.relations}, we get:

\begin{corollary}
\label{cor:compactproj}Let $\mathbb{D}$ be a composable class of filters.
Let $X$ be a \emph{topological} space. The following are equivalent:

\begin{enumerate}
\item X is $\mathbb{D}$-compact;

\item X is $\left( \int \mathbb{D}\right) $-compact;

\item for every $\left( \int \mathbb{D}\right) $-based convergence space $Y$%
, the projection $p_{Y}:X\times Y\rightarrow Y$ is $\left( \int \mathbb{D}%
\right) $-perfect;

\item for every $\mathbb{D}$-based atomic topological space $Y,$ the
projection $p_{Y}:X\times Y\rightarrow Y$ is closed.
\end{enumerate}
\end{corollary}

A similar result \cite[Theorem 1]{vaughan.proj} has been obtained by J.
Vaughan for topological spaces. He used nets instead of filters. To a class $%
\Omega $ of directed sets, we can associate a class $\mathbb{D}_{\Omega }$
of filters by 
\begin{equation*}
\mathcal{F}\in \mathbb{D}_{\Omega }\Longleftrightarrow \exists D\in \Omega
,\exists f:D\rightarrow \mathcal{F}:d\leq d^{\prime }\Longrightarrow
f(d^{\prime })\subset f(d).
\end{equation*}

The $\Omega $-\emph{net spaces} of \cite{vaughan.proj} are topological
spaces $(X,\xi )$ such that $\xi =T\mathrm{Base}_{\mathbb{D}%
_{\Omega }}\xi ;$ $\Omega $-\emph{Fr\'{e}chet spaces }are topological spaces 
$(X,\xi )$ such that $\xi =P\mathrm{Base}_{\mathbb{D}_{\Omega
}}\xi $ and $\Omega $-\emph{neighborhood spaces }are topological spaces $%
(X,\xi )$ such that $\xi =\mathrm{Base}_{\mathbb{D}_{\Omega }}\xi
. $ It is easy to show that a subspace of a topological space $(X,\xi )$
satisfying $\xi =T\mathrm{Base}_{\mathbb{D}}\xi $ is $\left( \int 
\mathbb{D}\right) $-based (see for instance \cite{DM.products}). Therefore 
\cite[Theorem 1]{vaughan.proj} follows from Corollary \ref{cor:compactproj}.
In particular, when $\mathbb{D}$ ranges over the classes $\mathbb{F}$ of all
filters, $\mathbb{F}_{\omega }$ of countably based filters, and $\mathbb{F}%
_{\wedge \omega }$ of countably deep filters (\footnote{%
A\ filter $\mathcal{F}$ is \emph{countably deep }if $\bigcap \mathcal{A\in F}
$ whenever $\mathcal{A}$ is a countable subfamily of $\mathcal{F}.$})
Corollary \ref{cor:compactproj} leads to

\begin{corollary}
(Mr\'{o}wka-Kuratowski \cite[Theorem 3.1.16]{Eng}) The following are
equivalent for a topological space $X$:

\begin{enumerate}
\item $X$ is compact;

\item $p_{Y}:X\times Y\rightarrow Y$ is perfect for every topological space $%
Y;$

\item $p_{Y}:X\times Y\rightarrow Y$ is closed for every topological space $%
Y.$
\end{enumerate}
\end{corollary}

\begin{corollary}
(Noble \cite[Corollary 2.4]{closedproj}) The following are equivalent for a
topological space $X$:

\begin{enumerate}
\item $X$ is countably compact;

\item $p_{Y}:X\times Y\rightarrow Y$ is countably perfect for every
subsequential (\footnote{%
A topological space is \emph{sequential} if every sequentially closed subset
is closed and \emph{subsequential} if it is homeomorphic to subspace of a
sequential space.}) topological space $Y;$

\item $p_{Y}:X\times Y\rightarrow Y$ is closed for every first-countable
topological space $Y.$
\end{enumerate}
\end{corollary}

\begin{corollary}
(Noble \cite[Corollary 2.3]{closedproj}) The following are equivalent for a
topological space $X$:

\begin{enumerate}
\item $X$ is Lindel\"{o}f;

\item $p_{Y}:X\times Y\rightarrow Y$ is inversely Lindel\"{o}f for every
topological $P$-space (\footnote{%
A topological space is \emph{a P-space } if every countable intersection of
open subsets is open; equivalently if it is $\mathbb{F}_{\wedge \omega }$%
-based.}) $Y;$

\item $p_{Y}:X\times Y\rightarrow Y$ is closed for every topological $P$%
-space $Y.$
\end{enumerate}
\end{corollary}

To a class $\mathbb{D}$ of filters, S. Dolecki associated in \cite{quest2}
two fundamental concrete functors of the category of convergence spaces: a
reflector $\mathrm{Adh}_{\mathbb{D}}$ where

\begin{equation}
\lm_{\mathrm{Adh}_{\mathbb{D}}\xi }\mathcal{F}=\bigcap_{%
\mathbb{D}\backepsilon \mathcal{D}\#\mathcal{F}}\mathrm{adh}_{\xi }%
\mathcal{D}  \label{eq:Adh}
\end{equation}

and a coreflector $\mathrm{Base}_{\mathbb{D}}$ where%
\begin{equation}
\lm_{\mathrm{Base}_{\mathbb{D}}\xi }\mathcal{F}=\bigcup_{%
\mathbb{D}\backepsilon \mathcal{D}\leq \mathcal{F}}\lm_{\xi }%
\mathcal{D}.  \label{eq:Base}
\end{equation}%
Applied to the case where $A$ is a singleton, Theorem \ref%
{th:Dcompactproduct} rephrases in convergence theoretic terms as follows.

\begin{theorem}
\label{th:AdhD}Let $\mathbb{D}$ be a composable class of filters and let $%
\xi $ and $\theta $ be two convergences on $X$. The following are equivalent:

\begin{enumerate}
\item $\theta \geq \mathrm{Adh}_{\mathbb{D}}\xi ;$

\item $\theta \times \mathrm{Base}_{\mathbb{D}}\tau \geq \mathrm{%
Adh}_{\mathbb{D}}\left( \xi \times \tau \right) $ for every convergence $%
\tau ;$

\item $\theta \times \tau \geq P\left( \xi \times \tau \right) $ for every $%
\mathbb{D}$-based atomic topology $\tau .$
\end{enumerate}
\end{theorem}

\begin{proof}
$(1\Longrightarrow 2).$ Let $x\in \lm_{\theta }\mathcal{F}$ and let $y\in
\lm_{\tau }\mathcal{G}$ with $\mathcal{G}\in \mathbb{D}.$ By assumption, $%
x\in \lm_{\mathrm{Adh}_{\mathbb{D}}\xi }\mathcal{F};$ in other words, $%
\mathcal{F}$ is $\mathbb{D}$-compact at $\{x\}$ and $\mathcal{G}\in \mathbb{D%
}$ is compact at $\{y\}$. By Theorem \ref{th:Dcompactproduct}, $\mathcal{F}%
\times \mathcal{G}$ is $\mathbb{D}$ -compact at $\{(x,y)\},$ that is, $%
(x,y)\in \lm_{\mathrm{Adh}_{\mathbb{D}}\left( \xi \times \tau \right) }(%
\mathcal{F}\times \mathcal{G}).$

$(2\Longrightarrow 3)$ is obvious and $(3\Longrightarrow 1)$ follows from $%
(3\Longrightarrow 1)$ in Theorem \ref{th:Dcompactproduct}. Indeed, if $x\in
\lm_{\theta }\mathcal{F},$then for every atomic topological space $(Y,\tau
) $ and every $\mathbb{D}$-filter $\mathcal{G}$ that converges to $y$ in $Y,$
$(x,y)\in \lm_{P\left( \xi \times \tau \right) }\left( \mathcal{F}\times 
\mathcal{G}\right) ,$ that is, the filter $\mathcal{F}\times \mathcal{G}$ is 
$\mathbb{F}_{1}$-compact at $\{(x,y)\}$, so that $\mathcal{F}$ is $\mathbb{D}
$-compact at $\{x\}.$ Hence $x\in \lm_{\mathrm{Adh}_{\mathbb{D}}\xi }%
\mathcal{F}.$
\end{proof}

The result above was essentially proved in \cite[Theorem 7.1]{mynard} but
was not stated explicitely in \cite{mynard}.

Let $\mathbb{D}$ and $\mathbb{J}$ be two classes of filters. A convergence
space is called $(\mathbb{J}/\mathbb{D})$-\emph{accessible} if whenever $x$
is an adherent point of a $\mathbb{J}$-filter $\mathcal{J},$ there exists a $%
\mathbb{D}$-filter $\mathcal{D}$ which converges to $x$ and meshes with $%
\mathcal{J}$. S. Dolecki introduced the notion (under a different name) in 
\cite{quest2} and noticed that when $\mathbb{D}$ is the class of countably
based filters and $\mathbb{J}$ ranges over the classes of all, of countably
deep, of countably based, of principal, of principal of closed sets filters,
then $(\mathbb{J}/\mathbb{D})$-accessible topological spaces are exactly the
bisequential \cite{mich.bi}, weakly bisequential \cite{liu.weaklybiseq},
strongly Fr\'{e}chet (countably bisequential in \cite{quest}), Fr\'{e}chet
and sequential spaces respectively. Additionally, he noticed that a
convegence $\xi $ is $(\mathbb{J}/\mathbb{D})$-accessible if and only if $%
\xi \geq \mathrm{Adh}_{\mathbb{J}}\mathrm{Base}_{\mathbb{D}}\xi $.

In view of \cite[Theorem 1]{myn.relations}, we obtain:

\begin{corollary}
\label{cor:JDaccess} Let $\mathbb{J\subset D}$ be two classes of filters
containing principal filters. Assume that a product of two $\mathbb{D}$%
-filters is a $\mathbb{D}$-filter. The following are equivalent:

\begin{enumerate}
\item $\xi $ is $\mathbb{(J}/\mathbb{D)}$-accessible;

\item $\xi \times \tau $ is $\mathbb{(J}/\mathbb{D)}$-accessible for every $%
\mathbb{J}$-based convergence space $(Y,\tau );$

\item $\xi \times \tau $ is $\mathbb{(F}_{1}/\mathbb{D)}$-accessible for
every atomic $\mathbb{J}$-based topological space $(Y,\tau ).$
\end{enumerate}
\end{corollary}

\begin{proof}
Notice that $\mathrm{Base}_{\mathbb{J}}\geq \mathrm{Base}%
_{\mathbb{D}}$ because $\mathbb{J\subset D}$.

$\left( 1\Longrightarrow 2\right) .$ If $\xi \geq \mathrm{Adh}_{%
\mathbb{J}}\mathrm{Base}_{\mathbb{D}}\xi $ and $\tau =\mathrm{Base%
}_{\mathbb{J}}\tau ,$ then 
\begin{equation*}
\xi \times \tau \geq \mathrm{Adh}_{\mathbb{J}}\mathrm{Base}%
_{\mathbb{D}}\xi \times \tau \geq \mathrm{Adh}_{\mathbb{J}%
}\mathrm{Base}_{\mathbb{D}}\xi \times \mathrm{Base}_{%
\mathbb{D}}\tau =\mathrm{Adh}_{\mathbb{J}}\mathrm{Base}_{%
\mathbb{D}}\left( \xi \times \tau \right) ,
\end{equation*}%
so that $\xi \times \tau $ is $\mathbb{(J}/\mathbb{D)}$-accessible.

$\left( 2\Longrightarrow 3\right) $ is clear because $\mathbb{F}_{1}\mathbb{%
\subset J}.$

$\left( 3\Longrightarrow 1\right) $ The convergence $\xi $ satisfies 
\begin{equation*}
\xi \times \tau \geq P\mathrm{Base}_{\mathbb{D}}\left( \xi \times
\tau \right) =P\left( \mathrm{Base}_{\mathbb{D}}\xi \times \tau
\right)
\end{equation*}%
for every $\mathbb{J}$-based atomic topology $\tau .$ By Theorem \ref%
{th:AdhD}, $\xi \geq \mathrm{Adh}_{\mathbb{J}}\mathrm{Base}%
_{\mathbb{D}}\xi .$
\end{proof}

In particular, when $\mathbb{J=D=F}_{\omega },$ it shows the following
generalization to convergence spaces of \cite[Propostions 4.D.4 and 4.D.5]%
{quest}:

\begin{corollary}
A convergence space is strongly Fr\'{e}chet if and only if its product with
every first-countable convergence (equivalently, every atomic
first-countable topological space) is strongly Fr\'{e}chet (equivalently Fr%
\'{e}chet).
\end{corollary}

An $\mathbb{F}_{1}$-based convergence is called \emph{finitely generated. }%
Finitely generated topologies are often called \emph{Alexandroff topologies. 
}When $\mathbb{J=F}_{1}$ and $\mathbb{D=F}_{\omega },$ Corollary \ref%
{cor:JDaccess} particularizes to

\begin{corollary}
\cite{mynard} A topological (or convergence) space is Fr\'{e}chet if and
only if its product with every finitely generated convergence space
(equivalently, Alexandroff topology) is Fr\'{e}chet .
\end{corollary}

On the other hand, applying Theorem \ref{th:Dcompactproduct} for the image
of a general filter under a relation, we obtain the following corollary for
(possibly multi-valued) maps.

\begin{corollary}
Let $\mathbb{D}$ be a composable class of filters and let $%
R:X\rightrightarrows Z$. The following are equivalent:

\begin{enumerate}
\item $R$ is a $\mathbb{D}$-compact relation;

\item $R\times Id_{Y}:X\times Y\rightrightarrows Z\times Y$ is a $\mathbb{D}$%
-compact relation for every $\mathbb{D}$-based convergence space $Y;$

\item $R\times Id_{Y}:X\times Y\rightrightarrows Z\times Y$ is an $\mathbb{F}%
_{1}$-compact relation for every atomic $\mathbb{D}$-based topological space 
$Y.$
\end{enumerate}
\end{corollary}

In view of \cite[Theorem 13]{myn.relations}, the last result leads to:

\begin{corollary}
Let $\mathbb{D}$ be a composable class of filters, let $X$ be a topological
space, and let $f:X\rightarrow Y$ be a surjective map. The following are
equivalent:

\begin{enumerate}
\item $f$ is $\mathbb{D}$-perfect;

\item $f\times Id_{W}$ is $\mathbb{D}$-perfect for every $\mathbb{D}$-based
convergence space $W;$

\item $f\times Id_{W}$ is $\left( \int \mathbb{D}\right) $-perfect for every 
$\left( \int \mathbb{D}\right) $-based topological space $W;$

\item $f\times Id_{W}$ is closed for every $\mathbb{D}$-based topological
space $W.$
\end{enumerate}
\end{corollary}

In particular, \cite[Corollary 3.5 (iii), (iv), (v) and (vi)]{closproj2} are
special cases:

\begin{corollary}
\label{cor:Dperfectchar} Let $X$ be a topological space, and let $%
f:X\rightarrow Y$ be a surjective map.The following are equivalent:

\begin{enumerate}
\item $f$ is perfect;

\item $f\times Id_{W}$ is perfect for every topological space $W;$

\item $f\times Id_{W}$ is closed for every topological space $W.$
\end{enumerate}
\end{corollary}

\begin{corollary}
Let $X$ be a topological space, and let $f:X\rightarrow Y$ be a surjective
map.The following are equivalent:

\begin{enumerate}
\item $f$ is countably perfect;

\item $f\times Id_{W}$ is countably perfect for every subsequential
topological space $W;$

\item $f\times Id_{W}$ is closed for every first-countable topological space 
$W.$
\end{enumerate}
\end{corollary}

\begin{corollary}
Let $X$ be a topological space, and let $f:X\rightarrow Y$ be a surjective
map.The following are equivalent:

\begin{enumerate}
\item $f$ is inversely Lindel\"{o}f;

\item $f\times Id_{W}$ is inversely Lindel\"{o}f for every topological $P$%
-space $W;$

\item $f\times Id_{W}$ is closed for every topological $P$-space $W.$
\end{enumerate}
\end{corollary}

Similarily, in view of of \cite[Theorem 14]{myn.relations} stating that a
map $f:(X,\xi )\rightarrow (Y,\tau )$ is $\mathbb{D}$-quotient if and only
if $f:(X,f^{-}\tau )\rightarrow (Y,f\xi )$ is $\mathbb{D}$-compact (%
\footnote{$f\xi $ denotes the final convergence on $Y$ associated to $%
f:(X,\xi )\rightarrow Y,$ that is, the finest convergence on $Y$ making $f$
continuous. Dually, $f^{-}\tau $ denotes the initial convergence on $X$
associated to $f:X\rightarrow (Y,\tau ),$ that is, the coarsest convergence
on $X$ making $f$ continuous.}), we obtain:

\begin{corollary}
\label{cor:Dquot} Let $\mathbb{D}$ be a composable class of filters and let $%
f:X\rightarrow Y$ be a surjective map. The following are equivalent:

\begin{enumerate}
\item $f$ is $\mathbb{D}$-quotient;

\item $f\times Id_{W}$ is $\mathbb{D}$-quotient for every $\mathbb{D}$-based
convergence space $W;$

\item $f\times Id_{W}$ is hereditarily quotient for every $\mathbb{D}$-based
topological space $W.$
\end{enumerate}
\end{corollary}

Notice that even if $X$ and $Y$ are topological, the final convergence $f\xi 
$ may not be. Therefore, $\mathbb{D}$-quotientness and $\left( \int \mathbb{D%
}\right) $-quotientness are not equivalent. Special instances include the
following:

\begin{corollary}
(Michael \cite{mich.bi}) The following are equivalent for a surjective map $%
f:X\rightarrow Y$ :

\begin{enumerate}
\item $f$ is biquotient;

\item $f\times Id_{W}$ is biquotient for every convergence space $W;$

\item $f\times Id_{W}$ is hereditarily quotient for every topological space $%
W.$
\end{enumerate}
\end{corollary}

\begin{corollary}
(Michael \cite[Propositions 4.3 and 4.4]{quest}) The following are
equivalent for a surjective map $f:X\rightarrow Y$ :

\begin{enumerate}
\item $f$ is countably biquotient;

\item $f\times Id_{W}$ is countably biquotient for every first-countable
convergence space $W;$

\item $f\times Id_{W}$ is hereditarily quotient for every first-countable
topological space $W.$
\end{enumerate}
\end{corollary}

\emph{Weakly biquotient maps }\cite{liu.weaklybiseq} coincide with $\mathbb{F%
}_{\wedge \omega }$-quotient maps so that when $\mathbb{D=F}_{\wedge \omega
} $ Corollary \ref{cor:Dquot} specializes to:

\begin{corollary}
The following are equivalent for a surjective map $f:X\rightarrow Y$ :

\begin{enumerate}
\item $f$ is weakly biquotient;

\item $f\times Id_{W}$ is weakly biquotient for every $\mathbb{F}_{\wedge
\omega }$-based convergence space $W;$

\item $f\times Id_{W}$ is hereditarily quotient for every topological $P$%
-space $W.$
\end{enumerate}
\end{corollary}

Finally, since a multivalued map $R:X\rightrightarrows Y$ between two
topological spaces is \emph{upper semicontinuous (usc)} if and only if it is
an $\mathbb{F}_{1}$-compact relation and \emph{compact-valued upper
semicontinuous (usco) }if and only if it is an $\mathbb{F}$-compact
relation, we have:

\begin{corollary}
Let $R:X\rightrightarrows Y$ be a multivalued map between two topological
space. Then

\begin{enumerate}
\item $R$ is an usco map if and only if $R\times Id_{W}:X\times
W\rightrightarrows Y\times W$ is a usc map (equivalently usco map) for every
topological space $W$;

\item $R$ is an usc map if and only if $R\times Id_{W}:X\times
W\rightrightarrows Y\times W$ is a usc map for every Alexandroff topological
space $W.$
\end{enumerate}
\end{corollary}

\section{Products of $\mathbb{D}$-compact filters}

In Section \ref{sec:Dcompprod}, $\mathbb{D}$-compact filters are
characterized as those filters whose product with every compact $\mathbb{D}$%
-filters is $\mathbb{D}$-compact. In this section, we consider the following
related question: What are the filters whose product with every $\mathbb{D}$%
-compact filter (of a given class $\mathbb{J)}$ is $\mathbb{D}$-compact ?

\subsection{Compactly meshable filters}

The question above was answered in \cite{JLM}, where a simplified version of
the following notion was introduced :

A filter $\mathcal{F}$ is $\mathbb{M}$-\emph{compactly} $\mathbb{J}$ \emph{to%
} $\mathbb{D}$ \emph{meshable }at $A,$ or $\mathcal{F}$ is an $\mathbb{M}$%
-compactly $\left( \mathbb{J}/\mathbb{D}\right) _{\#}$-filter at $A$, if 
\begin{equation*}
\mathcal{J}\in \mathbb{J},\mathcal{J}\#\mathcal{F}\Longrightarrow \exists 
\mathcal{D}\in \mathbb{D},\mathcal{D}\#\mathcal{J}\text{ and }\mathcal{D}%
\text{ is }\mathbb{M}\text{-compact at }A.
\end{equation*}

Before proceeding with applications, recall (see \cite{myn.relations} for
details) that the notion of an $\mathbb{M}$-compactly $\left( \mathbb{J}/%
\mathbb{D}\right) _{\#}$-filter is instrumental not only in answering the
question above but also in characterizing a large number of classical
concepts. It generalizes the notions of total countable compactness in the
sense of Z. Frol\'{\i}k \cite{frolikpseudo} and more generally of total $%
\mathbb{D}$-compactness in the sense of J. Vaughan \cite{vaughan78} from
sets to filters (\footnote{%
Let $\mathbb{D}$ be a class of filters. A topological space is \emph{totally 
}$\mathbb{D}$-\emph{compact }if every $\mathbb{D}$-filter has a finer
(relatively) compact $\mathbb{D}$-filter. It is easy to see that if $\mathbb{%
D}$ is stable under finite supremum, then $\{X\}$ is $\mathbb{F}$-compactly $%
\mathbb{D}$ to $\mathbb{D}$ meshable (at $X)$ iff $X$ is totally $\mathbb{D}$%
-compact. The notion of total countable compactness corresponds to $\mathbb{%
D=F}_{\omega }.$}).

\begin{proposition}
\label{pro:local M-compactoidly}\cite[Proposition 15]{myn.relations} Let $%
\mathbb{D}$, $\mathbb{J}$ and $\mathbb{M}$ be three classes of filters, and
let $\xi $ and $\theta $ be two convergences on $X$. The following are
equivalent:

\begin{enumerate}
\item $\theta \geq \mathrm{Adh}_{\mathbb{J}}\mathrm{Base}_{\mathbb{D}}%
\mathrm{Adh}_{\mathbb{M}}\xi ;$

\item $\mathcal{F}$ is an $\mathbb{M}$-compactly $(\mathbb{J}/\mathbb{D}%
)_{\#}$ filter at $\{x\}$ in $\xi $ whenever $x\in \lm_{\theta }\mathcal{F}%
. $

In particular, $\xi =\mathrm{Adh}_{\mathbb{M}}\xi $ is $(\mathbb{J}/\mathbb{%
D})$-accessible if and only if $\mathcal{F}$ is an $\mathbb{M}$-compactly $(%
\mathbb{J}/\mathbb{D})_{\#}$-filter at $\{x\}$ whenever $x\in \lm \mathcal{F%
}.$
\end{enumerate}
\end{proposition}

In view of \cite{quest2}, this means that Fr\'{e}chet, strongly Fr\'{e}chet,
productively Fr\'{e}chet, weakly bisequential, bisequential and radial
topological spaces among others, can be characterized in terms of $\mathbb{M}
$-compactly $(\mathbb{J}/\mathbb{D})_{\#}$-filters relative to a singleton,
for various instances of $\mathbb{J},$ $\mathbb{D}$ and $\mathbb{M}.$
Characterizations of Fr\'{e}chet and strongly Fr\'{e}chet spaces in terms
similar in spirit to those in Proposition \ref{pro:local M-compactoidly}
were obtained in \cite{active}. We take this opportunity to acknowledge that
even though productively Fr\'{e}chet spaces were not fully characterized in 
\cite{active}, important ideas and tools at work in \cite{JM} and \cite%
{mynard-jordan} were already introduced in \cite{active}.

More generally, the notion is instrumental in characterizing a number of
classes of maps. A \emph{relation} $R:(X,\xi )\rightrightarrows (Y,\tau )$ 
\emph{is} $\mathbb{M}$-\emph{compactly} $\left( \mathbb{J}/\mathbb{D}\right) 
$\emph{-meshable} if 
\begin{equation*}
\mathcal{F}\underset{\xi }{\rightarrow }x\Longrightarrow R(\mathcal{F})\text{
is }\mathbb{M}\text{-compactly }\left( \mathbb{J}/\mathbb{D}\right) \text{%
-meshable at }Rx\text{ in }\tau .
\end{equation*}

\begin{theorem}
\label{th:Mquot+range}\cite[Theorem 16]{myn.relations} Let $\mathbb{M\subset
J}$, let $\tau =\mathrm{Adh}_{\mathbb{M}}\tau $ and let $f:(X,\xi
)\rightarrow (Y,\tau )$ be a continuous surjection. The map $f$ is $\mathbb{M%
}$-quotient with $\mathbb{(J}/\mathbb{D)}$-accessible range if and only if $%
f:(X,f^{-}\tau )\rightarrow (Y,f\xi )$ is an $\mathbb{M}$-compactly $\left( 
\mathbb{J}/\mathbb{D}\right) $-meshable relation.
\end{theorem}

A convergence $\xi $ is $P$-\emph{diagonal}, if $\lm_{\xi }\mathcal{F}%
\subset \lm_{\xi }\int_{\mathcal{F}}\mathcal{V}_{\xi }(\mathcal{\cdot })$
for every filter $\mathcal{F}.$ The notation $\mathrm{adh}_{\xi
}^{\natural }(\mathbb{M)\subset M}$ means that the filter generated by $\{%
\mathrm{adh}_{\xi }M:M\in \mathcal{M\}}$ is in the class $\mathbb{M%
}$ whenever $\mathcal{M}$ is.

\begin{theorem}
\label{th:Mperfect+range}\cite[Theorem 17]{myn.relations} Let $\mathbb{%
M\subset J}$ and $\mathbb{D}$ be three classes of filters, where $\mathbb{J}$
and $\mathbb{D}$ are $\mathbb{F}_{1}$-composable. Let $\tau =\mathrm{Adh}%
_{\mathbb{M}}\tau $ and let $\xi $ be a $P$-diagonal convergence
such that $\mathrm{adh}_{\xi }^{\natural }(\mathbb{M)\subset M}$ .
Let $f:(X,\xi )\rightarrow (Y,\tau )$ be a continuous surjection. The map $f$
is $\mathbb{M}$-perfect with $\mathbb{(J}/\mathbb{D)}$-accessible range if
and only if $f^{-}:(Y,\tau )\rightrightarrows (X,\xi )$ is an $\mathbb{M}$%
-compactly $\left( \mathbb{J}/\mathbb{D}\right) $-meshable relation.
\end{theorem}

The following tables gather instances of these two results, for various
classes $\mathbb{M},$ $\mathbb{D}$ and $\mathbb{J}.$

\begin{center}
\begin{tabular}{|c|c|c|c|}
\hline
$\mathbb{M}$ & $\mathbb{J}$ & $\mathbb{D}$ & map $f$ as in Theorem \ref%
{th:Mquot+range} \\ \hline
$\mathbb{F}_{1}{}$ & $\mathbb{F}{}$ & $\mathbb{F}_{1}{}$ & hereditarily
quotient with finitely generated range \\ \hline
$\mathbb{F}_{1}{}$ & $\mathbb{F}_{1}{}$ & $\mathbb{F}_{\omega }{}$ & 
hereditarily quotient with Fr\'{e}chet range \\ \hline
$\mathbb{F}_{1}{}$ & $\mathbb{F}_{\omega }{}$ & $\mathbb{F}_{\omega }{}$ & 
hereditarily quotient with strongly Fr\'{e}chet range \\ \hline
$\mathbb{F}_{1}{}$ & $\mathbb{F}{}$ & $\mathbb{F}_{\omega }{}$ & 
hereditarily quotient with bisequential range \\ \hline
$\mathbb{F}_{1}$ & $\mathbb{F}$ & $\mathbb{F}$ & hereditarily quotient \\ 
\hline
$\mathbb{F}_{\omega }$ & $\mathbb{F}_{\omega }{}$ & $\mathbb{F}_{1}$ & 
countably biquotient with finitely generated range \\ \hline
$\mathbb{F}_{\omega }$ & $\mathbb{F}_{\omega }$ & $\mathbb{F}_{\omega }{}$ & 
countably biquotient with strongly Fr\'{e}chet range \\ \hline
$\mathbb{F}_{\omega }{}$ & $\mathbb{F}{}$ & $\mathbb{F}_{\omega }{}$ & 
countably biquotient with bisequential range \\ \hline
$\mathbb{F}_{\omega }{}$ & $\mathbb{F}{}$ & $\mathbb{F}{}$ & countably
biquotient \\ \hline
$\mathbb{F}{}$ & $\mathbb{F}{}$ & $\mathbb{F}_{1}{}$ & biquotient with
finitely generated range \\ \hline
$\mathbb{F}{}$ & $\mathbb{F}{}$ & $\mathbb{F}_{\omega }{}$ & biquotient with
bisequential range \\ \hline
$\mathbb{F}{}$ & $\mathbb{F}{}$ & $\mathbb{F}{}$ & biquotient \\ \hline
\end{tabular}

\begin{tabular}{|c|c|c|c|}
\hline
$\mathbb{M}$ & $\mathbb{J}$ & $\mathbb{D}$ & map $f$ as in Theorem \ref%
{th:Mperfect+range} \\ \hline
$\mathbb{F}_{1}{}$ & $\mathbb{F}{}$ & $\mathbb{F}_{1}{}$ & closed with
finitely generated range \\ \hline
$\mathbb{F}_{1}{}$ & $\mathbb{F}_{1}{}$ & $\mathbb{F}_{\omega }{}$ & closed
with Fr\'{e}chet range \\ \hline
$\mathbb{F}_{1}{}$ & $\mathbb{F}_{\omega }{}$ & $\mathbb{F}_{\omega }{}$ & 
closed with strongly Fr\'{e}chet range \\ \hline
$\mathbb{F}_{1}{}$ & $\mathbb{F}{}$ & $\mathbb{F}_{\omega }{}$ & closed with
bisequential range \\ \hline
$\mathbb{F}_{1}$ & $\mathbb{F}$ & $\mathbb{F}$ & closed \\ \hline
$\mathbb{F}_{\omega }$ & $\mathbb{F}_{\omega }{}$ & $\mathbb{F}_{1}$ & 
countably perfect with finitely generated range \\ \hline
$\mathbb{F}_{\omega }$ & $\mathbb{F}_{\omega }$ & $\mathbb{F}_{\omega }{}$ & 
countably perfect with strongly Fr\'{e}chet range \\ \hline
$\mathbb{F}_{\omega }{}$ & $\mathbb{F}{}$ & $\mathbb{F}_{\omega }{}$ & 
countably perfect with bisequential range \\ \hline
$\mathbb{F}_{\omega }{}$ & $\mathbb{F}{}$ & $\mathbb{F}{}$ & countably
perfect \\ \hline
$\mathbb{F}{}$ & $\mathbb{F}{}$ & $\mathbb{F}_{1}{}$ & perfect with finitely
generated range \\ \hline
$\mathbb{F}{}$ & $\mathbb{F}{}$ & $\mathbb{F}_{\omega }{}$ & perfect with
bisequential range \\ \hline
$\mathbb{F}{}$ & $\mathbb{F}{}$ & $\mathbb{F}{}$ & perfect \\ \hline
\end{tabular}
\end{center}

\subsection{Main Product Theorem}

The purpose of the remaining part of the paper is now to present
applications of Theorem \ref{th:main} below. It is formally more general and
also considerably more complicated than \cite[Theorem 1]{JLM} in order to
accomodate more applications. It is however based on the same idea.
Moreover, similar results were obtained jointly with F. Jordan (Georgia
Southern University) and I. Labuda (University of Mississippi) but were not
kept in full generality in \cite{JLM}. I thank them both for their
contributions to this subsection.

\begin{lemma}
\label{lem:compactprod}If $\mathcal{F}$ is a compact filter (at $\mathcal{A}%
) $ on $X$ such that $\mathcal{MF}\in \mathbb{M(}Y\mathbb{)}$ for every $%
\mathcal{M}\in \mathbb{M(}X\times Y\mathbb{)},$ and $\mathcal{G}$ is $%
\mathbb{M}$-compact (at $\mathcal{B}),$ then $\mathcal{F}\times \mathcal{G}$
is $\mathbb{M}$-compact (at $\mathcal{A\times B}).$
\end{lemma}

\begin{proof}
Let $\mathcal{M}$ be an $\mathbb{M}$-filter such that $\mathcal{M}\#\mathcal{%
F}\times \mathcal{G}.$ Then $\mathcal{M}\left( \mathcal{F}\right) \#\mathcal{%
G}$ and $\mathcal{MF}$ is an $\mathbb{M}$-filter, so that there exists $%
\mathcal{U}\#\mathcal{M}\left( \mathcal{F}\right) $ such that $\mathcal{U}%
\rightarrow y.$ The filter $\mathcal{M}^{-}\left( \mathcal{U}\right) $
meshes with the compact filter $\mathcal{F}$ and so there exists $\mathcal{W}%
\#\mathcal{M}^{-}\left( \mathcal{U}\right) $ such that $\mathcal{%
W\rightarrow }x.$ Then $(x,y)\in \mathrm{adh}\mathcal{M}.$
\end{proof}

In particular, if $\mathcal{F}=\{X\}$ and $\mathcal{G}=\{Y\},$ it shows that
the product of a compact space with a countably compact (resp. Lindel\"{o}f,
pseudocompact) space is countably compact (resp. Lindel\"{o}f,
pseudocompact).

\begin{theorem}
\label{th:main} Let $\mathbb{D}$ and $\mathbb{M}$ be two composable classes
of filters containing principal filters and let $\mathbb{J}$ and $\mathbb{K}$
be two $\mathbb{D}$-composable classes of filters. Let $\mathcal{F}\in 
\mathbb{K(}X)$ and $A\subset X.$ The following are equivalent:

\begin{enumerate}
\item $\mathcal{F}$ is a $\mathbb{M}$-compactly $(\mathbb{J}/\mathbb{D}%
)_{\#} $- filter at $A\subset X;$

\item for every $Y$, every $B\subset Y$ and every $\left( \mathbb{K}/\mathbb{%
J}\right) _{\#\geq }$-filter $\mathcal{G}$ which is a compactly $(\mathbb{D}/%
\mathbb{M})_{\#}$-filter at $B,$ the filter $\mathcal{F}\times \mathcal{G}$
is an $\mathbb{M}$-compactly $(\mathbb{D}/\mathbb{D}\times \mathbb{M})_{\#}$%
-filter at $A\times B;$

\item for every $(\mathbb{D}/\mathbb{M})$-accessible space $Y$, every $%
B\subset Y$ and every $\mathbb{J}$-filter $\mathcal{G}$ which is $\mathbb{D}$%
-compact at $B,$ the filter $\mathcal{F}\times \mathcal{G}$ is $\left( 
\mathbb{D\cap M}\right) $-compact at $A\times B;$

\item for every $\mathbb{M}$-based convergence space $Y$ and $y\in Y,$ and
for every $\mathbb{J}$-filter $\mathcal{G}$ which is $\mathbb{D}$-compact at 
$\{y\},$ the filter $\mathcal{F}\times \mathcal{G}$ is $\mathbb{F}_{1}$%
-compact at $A\times \{y\};$

\item for every $\mathbb{M}$-based (possibly non-Hausdorff) topological
space $Y$ and $B\subset Y,$ and for every $\mathbb{J}$-filter $\mathcal{G}$
which is $\mathbb{D}$-compact at $B,$ the filter $\mathcal{F}\times \mathcal{%
G}$ is $\mathbb{F}_{1}$-compact at $A\times B.$\newline
\end{enumerate}
\end{theorem}

\begin{proof}
$\left( 1\Longrightarrow 2\right) $ Let $\mathcal{D}$ be a $\mathbb{D}$%
-filter such that $\mathcal{D}\#\mathcal{F}\times \mathcal{G}.$ We can
assume without loss of generality that $\mathcal{G}\in \mathbb{J}.$ Indeed, $%
\mathcal{D}\left( \mathcal{F}\right) \in \mathbb{K}$ because $\mathbb{K}$ is 
$\mathbb{D}$-composable, and $\mathcal{D}\left( \mathcal{F}\right) \#%
\mathcal{G}$. Therefore, there exists a $\mathbb{J}$-filter $\mathcal{G}%
^{\prime }\#\mathcal{D}\left( \mathcal{F}\right) $ and finer than $\mathcal{G%
}.$ Since $\mathcal{G}^{\prime }\geq \mathcal{G},$ the filter $\mathcal{G}%
^{\prime }$ is a compactly $(\mathbb{D}/\mathbb{M})_{\#}$-filter at $B.$
Moreover, $\mathcal{D}\#\left( \mathcal{F}\times \mathcal{G}^{\prime
}\right) .$ From now on, assume that $\mathcal{G}\in \mathbb{J}.$ As $%
\mathbb{J}$ is $\mathbb{D}$-composable, $\mathcal{D}^{-}\left( \mathcal{G}%
\right) $ is a $\mathbb{J}$-filter and $\mathcal{D}^{-}\left( \mathcal{G}%
\right) \#\mathcal{F}.$ Since $\mathcal{F}$ is an $\mathbb{M}$-compactly $(%
\mathbb{J}/\mathbb{D})_{\#}$filter at $A$, there exists $\mathbb{D}$-filter $%
\mathcal{C}\#\mathcal{D}^{-}\left( \mathcal{G}\right) $ which is $\mathbb{M}$%
-compact at $A.$ Now $\mathcal{D}\left( \mathcal{C}\right) \#\mathcal{G}$
and $\mathcal{D}\left( \mathcal{C}\right) $ is a $\mathbb{D}$-filter, so
that there exists a filter $\mathcal{M}$ in $\mathbb{M}$ which is compact at 
$B$ and meshes with $\mathcal{D}\left( \mathcal{C}\right) $. By Lemma \ref%
{lem:compactprod}, $\mathcal{C}\times \mathcal{M}$ is an $\mathbb{M}$%
-compact filter at $A\times B$ meshing with $\mathcal{D}$ because $\mathbb{M}
$ is composable. Moreover, $\mathcal{C}\times \mathcal{M}\in \mathbb{D\times
M}$. Hence, $\mathcal{F}\times \mathcal{G}$ is an $\mathbb{M}$-compactly $(%
\mathbb{D}/\mathbb{D\times M})_{\#}$-filter at $A\times B$.

$\left( 2\Longrightarrow 3\right) $ because a $\mathbb{D}$-compact filter in
a $\mathbb{(D}/\mathbb{M)}$-accessible space is compactly $\mathbb{(D}/%
\mathbb{M)}$ meshable and a $\mathbb{M}$-compactly $(\mathbb{D}/\mathbb{%
D\times M})_{\#}$-filter is also $\left( \mathbb{D\cap M}\right) $-compact.

$\left( 3\Longrightarrow 4\right) $ and $\left( 3\Longrightarrow 5\right) $
are clear, as $\mathbb{F}_{1}\subset \mathbb{M\cap D}$ and every $\mathbb{M}$%
-based convergence space is $(\mathbb{D}/\mathbb{M})$-accessible.

$\left( 4\Longrightarrow 1\right) $. If $\mathcal{F}$ is not $\mathbb{M}$%
-compactly $(\mathbb{J}/\mathbb{D})_{\#}$ at $A,$ then there exists a $%
\mathbb{J}$-filter $\mathcal{J}\#\mathcal{F}$ such that for every $\mathbb{D}
$-filter $\mathcal{D}\#\mathcal{J},$ there exists a $\mathbb{M}$-filter $%
\mathcal{M}_{\mathcal{D}}\#\mathcal{D}$ such that $\mathrm{adh}\mathcal{M}_{%
\mathcal{D}}\cap A=\emptyset .$ Pick $y_{0}\in A$ and denote by $Y$ a copy
of $X$ endowed with the atomic $\mathbb{M}$-based convergence structure
defined by $y_{0}\in \lm \mathcal{G}$ iff there exists $\mathcal{D}\#%
\mathcal{J}$ such that $\mathcal{G}\geq \mathcal{M}_{\mathcal{D}}\wedge
\{y_{0}\}.$ Then $\mathcal{J}$ is $\mathbb{D}$-compact at $\{y_{0}\}$ in $Y$%
, but $\mathcal{F}\times \mathcal{J}$ is not $\mathbb{F}_{1}$-compact at $%
A\times \{y_{0}\}.$ Indeed, $\Delta =\{(x,x):x\in X,x\neq y_{0}\}\subset
X\times Y$ is in $\mathbb{F}_{1}$ and $\Delta \#\left( \mathcal{F}\times 
\mathcal{J}\right) $ because $\mathcal{F}\#\mathcal{J}.$ But $\mathrm{adh}%
\Delta \cap A\times \{y_{0}\}=\varnothing .$ Indeed, a filter on $\Delta $
can be assumed to be of the form $\mathcal{H}\times \mathcal{H}.$ Now if $%
\mathcal{H}$ converges to $\{y_{0}\}$ in $Y,$ then $\mathcal{H\geq M}_{%
\mathcal{D}}$ so that $\mathcal{H}$ cannot converge to $y_{0}\in A$ in $X,$
since $\mathrm{adh}\mathcal{M}_{\mathcal{D}}\cap A=\emptyset .$

$\left( 5\Longrightarrow 1\right) $. In the argument $\left(
4\Longrightarrow 1\right) ,$consider instead of the convergence space $Y,$
the $\mathbb{M}$-based topological space $Y$=$X\oplus \{\mathcal{M}_{%
\mathcal{D}}:\mathcal{D}\#\mathcal{J},$ $\mathcal{D}\in \mathbb{D}\}.$ Then
the filter $\widehat{\mathcal{J}}$ generated by $\mathcal{J}$ on $Y$ is $%
\mathbb{D}$-compact at $B=\{\mathcal{M}_{\mathcal{D}}:\mathcal{D}\#\mathcal{J%
},$ $\mathcal{D}\in \mathbb{D}\}$ but $\mathcal{F}\times \widehat{\mathcal{J}%
}$ is not $\mathbb{F}_{1}$-compact at $A\times B.$ Indeed, $\Delta
=\{(x,x):x\in X\}\subset X\times Y$ is in $\mathbb{F}_{1}$ and $\Delta
\#\left( \mathcal{F}\times \widehat{\mathcal{J}}\right) $ because $\mathcal{F%
}\#\mathcal{J}.$ But $\mathrm{adh}\Delta \cap A\times B=\varnothing .$
Indeed, a filter on $\Delta $ is of the form $\mathcal{H}\times \mathcal{H}.$
Now if $\mathcal{H}$ converges to some point $\{\mathcal{M}_{\mathcal{D}}\}$
in $Y,$ then $\mathcal{H\geq M}_{\mathcal{D}}$ and $\mathcal{H}$ cannot
converge to any point of $A$ in $X,$ since $\mathrm{adh}\mathcal{M}_{%
\mathcal{D}}\cap A=\emptyset .$
\end{proof}

From the viewpoint of convergence, there is no reason to distinguish between
a sequence and the filter generated by the family of its tails. Therefore,
in this paper, sequences are identified to their associated filter and we
will freely treat sequences as filters. For instance, given a filter $%
\mathcal{M},$ we consider the set $\mathcal{E(M)}=\{(x_{n})_{n\in \mathbb{N}%
}:(x_{n})_{n\in \mathbb{N}}\geq \mathcal{M}\}$ of free sequences finer than $%
\mathcal{M}$ by applying this convention.

\begin{lemma}
\label{lem:fanlike} Let $\mathcal{M}$ be a filter on $X.$ The filter $%
\mathcal{M}$ admits a finer free sequence ($\mathcal{E(M)}\neq \emptyset $)
if and only if for every family $(\mathcal{G}_{\alpha })_{\alpha \in I}$ of
free filters on $X$ such that $\mathcal{M}\geq \bigwedge_{\alpha \in
I}\mathcal{G}_{\alpha }$ there exists $\alpha _{0}\in I$ and $(x_{n})_{n\in 
\mathbb{N}}\geq \mathcal{M}$ such that $(x_{n})_{n\in \mathbb{N}}\#\mathcal{G%
}_{\alpha _{0}}.$ In particular, $\mathcal{M}\#\mathcal{G}_{\alpha _{0}}.$
\end{lemma}

\begin{proof}
Assume that there exists $(x_{n})_{n\in \mathbb{N}}\geq \mathcal{M}$ and
that for every $\alpha \in I,$ there exists $G_{\alpha }\in \mathcal{G}%
_{\alpha }$ such that $G_{\alpha }\notin \left( (x_{n})_{n\in \mathbb{N}%
}\right) ^{\#}.$ Since $\mathcal{G}_{\alpha }$ is free, there exists $%
G_{\alpha }^{\prime }\in \mathcal{G}_{\alpha }$ such that $G_{\alpha
}^{\prime }\cap \{x_{n}:n\in \mathbb{N}\}=\emptyset .$ Then $%
\bigcup_{\alpha \in I}G_{\alpha }^{\prime }\in
\bigwedge_{\alpha \in I}\mathcal{G}_{\alpha }$ but $%
\bigcup_{\alpha \in I}G_{\alpha }^{\prime }\notin (x_{n})_{n\in 
\mathbb{N}}.$ Therefore $(x_{n})_{n\in \mathbb{N}}\ngeq
\bigwedge_{\alpha \in I}\mathcal{G}_{\alpha }.$

The converse is obvious.
\end{proof}

We can now give an alternative version of $\left( 2\Longrightarrow 1\right) $
in Theorem \ref{th:main}.

\begin{proposition}
\label{pro:altconverse}Let $\mathbb{M}$ be a class of filters such that $%
\mathcal{E(M)}\neq \emptyset $ whenever $\mathcal{M}\in \mathbb{M}.$

Assume that for every $(\mathbb{F}_{1}/\mathbb{M)}$-accessible atomic
topological space $Y$ and every $\mathbb{J}$-filter $\mathcal{J},$ which is
compactly $\mathbb{D}$ to $\mathbb{M}$ meshable at the non-isolated point $%
\{\infty \}$ of $Y,$ the filter $\mathcal{F}\times \mathcal{J}$ is an $%
\mathbb{F}_{1}$-compactly $(\mathbb{F}_{1}/\mathbb{M)}_{\#}$-filter at $%
A\times \{\infty \}.$ Then $\mathcal{F}$ is an $\mathbb{M}$-compactly $(%
\mathbb{J}/\mathbb{D})_{\#}$- filter at $A.$
\end{proposition}

\begin{proof}
If $\mathcal{F}$ is not $\mathbb{M}$-compactly $(\mathbb{J}/\mathbb{D})_{\#}$
at $A,$ then there exists a $\mathbb{J}$-filter $\mathcal{J}\#\mathcal{F}$
such that for every $\mathbb{D}$-filter $\mathcal{D}\#\mathcal{J},$ there
exists an $\mathbb{M}$-filter $\mathcal{M}_{\mathcal{D}}\#\mathcal{D}$ such
that $\mathrm{adh}\mathcal{M}_{\mathcal{D}}\cap A=\emptyset .$ Let $%
Y=X\oplus \bigwedge \left\{ \mathcal{M}_{\mathcal{D}}:\mathcal{D}\in \mathbb{%
D},\mathcal{D}\#\mathcal{J}\right\} $ and denote by $\infty $ the point $%
\bigwedge \left\{ \mathcal{M}_{\mathcal{D}}:\mathcal{D}\in \mathbb{D},%
\mathcal{D}\#\mathcal{J}\right\} $ of $Y.$ Since infima of $\mathbb{M}$%
-filters are exactly $\left( \mathbb{F}_{1}/\mathbb{M}\right) _{\#\geq }$%
-filters, $Y$ is an $(\mathbb{F}_{1}/\mathbb{M)}$-accessible topological
space $Y.$ By definition, $\mathcal{J}$ is a compactly $(\mathbb{D}/\mathbb{M%
})_{\#}$-filter at $\{\infty \},$ but $\mathcal{F}\times \mathcal{J}$ is not
an $\mathbb{F}_{1}$-compactly $(\mathbb{F}_{1}/\mathbb{M)}_{\#}$-filter at $%
A\times \{\infty \}.$ Indeed, $\Delta =\{(x,x):x\in X\}$ meshes with $%
\mathcal{F}\times \mathcal{J}$ because $\mathcal{J}\#\mathcal{F}$ in $X.$ An 
$\mathbb{M}$-filter on $\Delta $ is of the form $\mathcal{M}\times \mathcal{M%
}$ where $\mathcal{M}\in \mathbb{M}(X).$ Assume that $\mathcal{M}\times 
\mathcal{M}\rightarrow (x,\infty )$ in $X\times Y.$ Then $\mathcal{M}\geq
\bigwedge_{\substack{ \mathcal{D}\in \mathbb{D}  \\ \mathcal{D}\#%
\mathcal{J}}}\mathcal{M}_{\mathcal{D}}$ because $\mathcal{M}\rightarrow
\infty .$ By Lemma \ref{lem:fanlike}, there exists a $\mathbb{D}$-filter $%
\mathcal{D}\#\mathcal{J}$ such that $\mathcal{M}\#\mathcal{M}_{\mathcal{D}}.$
Consequently, $x\notin A$ because $\mathrm{adh}_{X}\mathcal{M}_{\mathcal{D}%
}\cap A=\emptyset .$
\end{proof}

\section{Further applications}

\subsection{Global properties}

As observed in \cite{JLM}, the part $\left( 1\Longrightarrow 3\right) $ of
Theorem \ref{th:main} applied to principal filters $\mathcal{F}=\{X\}$ and $%
\mathcal{G}=\{Y\},$ for various instances of $\mathbb{D=J}$ and of $\mathbb{M%
}$ allows to recover results of J. Vaughan \cite{vaughan78}, and also to
provide new variants. For instance:

\begin{theorem}
\cite{JLM}

\begin{enumerate}
\item The product of a countably compact space and a compactly $\mathbb{(F}%
_{\omega }/\mathbb{F}_{\mathbb{\omega }})$-meshable space is countably
compact.

\item The product of a strongly Fr\'{e}chet countably compact space and a $%
\mathbb{F}_{\mathbb{\omega }}$-compactly $\mathbb{(F}_{\omega }/\mathbb{F}_{%
\mathbb{\omega }})$-meshable space is countably compact.
\end{enumerate}
\end{theorem}

For example, compact, sequentially compact, countably compact $k$-spaces are
all examples of compactly $\mathbb{(F}_{\omega }/\mathbb{F}_{\mathbb{\omega }%
})$-meshable space and every countably compact space is a $\mathbb{F}_{%
\mathbb{\omega }}$-compactly $\mathbb{(F}_{\omega }/\mathbb{F}_{\mathbb{%
\omega }})$-meshable space.

If $\mathcal{A}$ is a family of subsets of a convergence space $(X,\xi ),$
denote by $\mathcal{O}_{\xi }(\mathcal{A})$ the family $\{O$ open: $\exists
A\in \mathcal{A},$ $A\subset O\}.$ Accordingly, $\mathcal{O}_{\xi }\mathcal{(%
}\mathbb{D}\mathcal{)}$ will denote the class of $\mathbb{D}$-filters $%
\mathcal{D}$ such that $\mathcal{D}=\left( \mathcal{O}_{\xi }\mathcal{(D)}%
\right) ^{\uparrow }.$ A topological space $X$ is \emph{feebly compac}t if
and only if $\{X\}$ is $\mathbb{\mathcal{O}}(\mathbb{F}_{\omega })$-compact (%
\footnote{%
A Tychonoff space is feebly compact if and only if it is pseudocompact.}).

\begin{theorem}
\cite{JLM}

\begin{enumerate}
\item The product of a feebly compact space and a compactly $\mathbb{(%
\mathcal{O}}(\mathbb{F}_{\omega })/\mathcal{O}(\mathbb{F}_{\mathbb{\omega }%
}))$-meshable space is feebly compact.

\item The product of a $\mathbb{(\mathcal{O}}(\mathbb{F}_{\omega })/\mathbb{F%
}_{\mathbb{\omega }})$-accessible (in particular strongly Fr\'{e}chet$)$
feebly compact space and a $\mathbb{F}_{\mathbb{\omega }}$-compactly $%
\mathbb{(\mathcal{O}}(\mathbb{F}_{\omega })/\mathcal{O}(\mathbb{F}_{\mathbb{%
\omega }}))$-meshable space is feebly compact.
\end{enumerate}
\end{theorem}

\begin{theorem}
\cite{JLM}

\begin{enumerate}
\item The product of a Lindel\"{o}f space and a compactly $\mathbb{(F}%
_{\wedge \omega }/\mathbb{F}_{\wedge \mathbb{\omega }})$-meshable space is
Lindel\"{o}f.

\item The product of a weakly bisequential Lindel\"{o}f space and a $\mathbb{%
F}_{\mathbb{\omega }}$-compactly $\mathbb{(F}_{\wedge \omega }/\mathbb{F}%
_{\wedge \mathbb{\omega }})$-meshable space is Lindel\"{o}f.
\end{enumerate}
\end{theorem}

\subsection{ Local properties}

Theorem \ref{th:main} and Proposition \ref{pro:altconverse} applied in the
case of compactness at a singleton leads to the following.

\begin{theorem}
Let $\mathbb{D}\subset \mathbb{M}$ be two composable classes of filters
containing principal filters and assume that there exists a sequence $%
(x_{n})_{n\in \mathbb{N}}\geq \mathcal{M}$ whenever $\mathcal{M}\in \mathbb{M%
}.$ The following are equivalent for a topological space $X:$

\begin{enumerate}
\item $X$ is $\left( (\mathbb{D}/\mathbb{M})_{\#\geq }/\mathbb{D}\right) $%
-accessible;

\item $X\times Y$ is $(\mathbb{D}/\mathbb{M})$-accessible for every $(%
\mathbb{D}/\mathbb{M})$-accessible topological space $Y;$

\item $X\times Y$ is $(\mathbb{F}_{1}/\mathbb{M})$-accessible for every $(%
\mathbb{D}/\mathbb{M})$-accessible atomic topological space $Y.$
\end{enumerate}
\end{theorem}

\begin{proof}
$\left( 1\Longrightarrow 2\right) .$ Let $x\in \lm_{X}\mathcal{F}$ and $%
y\in \lm_{Y}\mathcal{G}.$ In view of Proposition \ref{pro:local
M-compactoidly}, $\mathcal{F}$ is an $\mathbb{M}$-compactly $\left( (\mathbb{%
D/M})_{\#\geq }/\mathbb{D}\right) _{\#}$-filter at $\{x\}$ and $\mathcal{G}$
is a compactly $(\mathbb{D/M})_{\#}$-filter at $\{y\}$ because $X$ is $%
\left( (\mathbb{D/M})_{\#\geq }/\mathbb{D}\right) $-accessible and $Y$ is $(%
\mathbb{D/M})$-accessible. Moreover, $\mathcal{G}$ can be assumed to be in $(%
\mathbb{D/M})_{\#\geq },$ which is a $\mathbb{D}$-composable class of
filters \cite{JM2}. By $(1\Longrightarrow 2)$ of Theorem \ref{th:main} with $%
\mathbb{J}=(\mathbb{D/M})_{\#\geq }$, $\mathcal{F}\times \mathcal{G}$ is an $%
\mathbb{M}$-compactly $\mathbb{(D/M)}_{\#}$-filter at$\{(x,y)\}$. Hence, $%
X\times Y$ is $(\mathbb{D/M})$-accessible.

$\left( 2\Longrightarrow 3\right) $ is trivial.

$\left( 3\Longrightarrow 1\right) .$ Let $x\in \lm_{X}\mathcal{F}$. We use
Proposition \ref{pro:altconverse} (with $\mathbb{J}=(\mathbb{D/M})_{\#\geq
}) $to show that $\mathcal{F}$ is an $\mathbb{M}$-compactly $\left( (\mathbb{%
D/M})_{\#\geq }/\mathbb{D}\right) _{\#}$-filter at $\{x\},$ which will show
that $X$ is $\left( (\mathbb{D/M})_{\#\geq }/\mathbb{D}\right) $-accessible
by Proposition \ref{pro:local M-compactoidly}. To this end, consider a $(%
\mathbb{D/M})_{\#\geq }$-filter $\mathcal{J}$ which is a compactly $(\mathbb{%
D/M})_{\#}$-filter at $\{\infty \}$ where $\infty $ is the non-isolated
point of an $(\mathbb{F}_{1}/\mathbb{M)}$-accessible atomic topological
space $Y.$ Notice that $\mathcal{J}$ is $\mathbb{F}_{1}$-compact at $%
\{\infty \},$ hence converges to $\infty $ in $Y$. Let $Y^{\prime }$ be the
(finer) atomic topological space obtained from $Y$ by letting $\mathcal{N}%
_{Y^{\prime }}(\infty )=\mathcal{J}\wedge \{\infty \}^{\uparrow }.$ The
space $Y^{\prime }$ is an atomic $(\mathbb{D/M})$-accessible topological
space, so that $X\times Y^{\prime }$ is an $(\mathbb{F}_{1}\mathbb{/M})$%
-accessible topological space. Therefore $\mathcal{F}\times \mathcal{J}$ is
an $\mathbb{F}_{1}$-compactly $(\mathbb{F}_{1}/\mathbb{M)}_{\#}$-filter at $%
\{(x,\infty )\}.$ By Proposition \ref{pro:altconverse}, $\mathcal{F}$ is an $%
\mathbb{M}$-compactly $\left( (\mathbb{D/M})_{\#\geq }/\mathbb{D}\right)
_{\#}$-filter at $\{x\}.$
\end{proof}

In particular if $\mathbb{D=M=F}_{\omega },$ we obtain:

\begin{corollary}
\cite{JM} The following are equivalent:

\begin{enumerate}
\item $X$ is productively Fr\'{e}chet;

\item $X\times Y$ is strongly Fr\'{e}chet for every strongly Fr\'{e}chet
topological space $Y;$

\item $X\times Y$ is Fr\'{e}chet for every atomic strongly Fr\'{e}chet
topological space $Y.$
\end{enumerate}
\end{corollary}

For $\mathbb{D=F}_{1}$ and $\mathbb{M=F}_{\omega },$ we obtain:

\begin{corollary}
\cite{mynard} The following are equivalent:

\begin{enumerate}
\item $X$ is \emph{finitely generated} (i.e., every point has a minimal
neighborhood);

\item $X\times Y$ is Fr\'{e}chet for every Fr\'{e}chet topological space $Y; 
$

\item $X\times Y$ is Fr\'{e}chet for every atomic Fr\'{e}chet topological
space $Y.$
\end{enumerate}
\end{corollary}

\subsection{Products of Maps}

In view of Theorems \ref{th:Mperfect+range} and \ref{th:Mquot+range} and of 
\cite[Theorems 13 and 14]{myn.relations}, Theorem \ref{th:main} has
important consequences in terms of product of maps. More specifically:

\begin{theorem}
\label{th:compactoidrelation}Let $\mathbb{D}$ and $\mathbb{M}$ be two
composable classes of filters containing principal filters and let $\mathbb{J%
}$ be a $\mathbb{D}$-composable class of filters. The following are
equivalent for a relation $R:X\rightrightarrows Y$ :

\begin{enumerate}
\item $R$ is an $\mathbb{M}$-compactly $(\mathbb{J}/\mathbb{D})$-meshable
relation$;$

\item for every compactly $\mathbb{(D}/\mathbb{M)}$-meshable relation $%
G:Z\rightrightarrows W$ where the convergence space $Z$ is $\mathbb{J}$%
-based, the relation $R\times G:X\times Z\rightrightarrows Y\times W$ is an $%
\mathbb{M}$-compactly $\mathbb{(D}/\mathbb{D\times M)}$-meshable relation;

\item for every $\mathbb{D}$-compact relation $G:Z\rightrightarrows W$ where
the convergence spaces $W$ and $Z$ are respectively $(\mathbb{D}/\mathbb{M})$%
-accessible and $\mathbb{J}$-based, the relation $R\times G:X\times
Z\rightrightarrows Y\times W$ is $\left( \mathbb{D\cap M}\right) $-compact ;

\item for every $\mathbb{M}$-based convergence space $W$, the relation $%
R\times Id:X\times \mathrm{Base}_{\mathbb{J}}\mathrm{Adh}_{\mathbb{D}%
}W\rightrightarrows Y\times W$ is $\mathbb{F}_{1}$-compact ;

\item for every map $g:W\rightarrow Z,$where $W$ is an $\mathbb{M}$-based
topological space and $Z$ is a $\mathbb{J}$-based atomic topological space,
whose inverse relation $g^{-}:Z\rightrightarrows W$ is $\mathbb{D}$-compact
, the relation $R\times g^{-}:X\times Z\rightrightarrows Y\times W$ is $%
\mathbb{F}_{1}$-compact .
\end{enumerate}
\end{theorem}

\begin{proof}
$\left( 1\Longrightarrow 2\right) $ Let $x\in \lm_{X}\mathcal{F}$ and $z\in
\lm_{Z}\mathcal{G}.$ We can assume $\mathcal{G}$ and hence $G(\mathcal{G)}$
to be $\mathbb{J}$-filters. By assumption, $G(\mathcal{G})$ is a $\mathbb{J}$%
-filter that is compactly $\mathbb{(D}/\mathbb{M)}$-meshable at $Gy,$ and $R(%
\mathcal{F})$ is $\mathbb{M}$-compactly $(\mathbb{J}/\mathbb{D})$-meshable
at $Rx.$ By Theorem \ref{th:main} $\left( 1\Longrightarrow 2\right) $, $R(%
\mathcal{F})\times G(\mathcal{G})$ is $\mathbb{M}$-compactly $\mathbb{(D}/%
\mathbb{D\times M)}$-meshable at $Rx\times Gy$ in $Y\times W.$

$(2\Longrightarrow 3)$, $(3\Longrightarrow 4)$ and $(3\Longrightarrow 5)$
are obvious.

$(4\Longrightarrow 1).$ In view of Theorem \ref{th:main}, it is sufficient
to show that $x\in \lm \mathcal{F}$ implies that $R\left( \mathcal{F}%
\right) \times \mathcal{G}$ is $\mathbb{F}_{1}$-compact at $Rx\times \{w\}$
in $Y\times W$ whenever $\mathcal{G}$ is a $\mathbb{D}$-compact at $\{w\}$ $%
\mathbb{J}$-filter, where $W$ is an $\mathbb{M}$-based convergence space.
Notice that $w\in \lm_{\mathrm{Base}_{\mathbb{J}}\mathrm{Adh}_{\mathbb{D}%
}W}\mathcal{G}$. Therefore $\left( R\times Id\right) \left( \mathcal{F}%
\times \mathcal{G}\right) =R\left( \mathcal{F}\right) \times \mathcal{G}$ is 
$\mathbb{F}_{1}$-compact at $Rx\times \{w\}$ in $Y\times W$ and the
conclusion follows.

$(5\Longrightarrow 1).$ In view of Theorem \ref{th:main}, it is sufficient
to show that $x\in \lm \mathcal{F}$ implies that $R\left( \mathcal{F}%
\right) \times \mathcal{G}$ is $\mathbb{F}_{1}$-compact at $Rx\times B$
whenever $\mathcal{G}$ is $\mathbb{D}$-compact at $B\subset W,$ where $W$ is
an $\mathbb{M}$-based topological space. For each such $\mathcal{G},$
consider the relation $G_{\mathcal{G}}:Z\rightrightarrows W,$ where $Z$ is
the $\mathbb{J}$-based atomic topological space $W\oplus \{\mathcal{G}\}$,
defined by $G_{\mathcal{G}}(w)=\{w\}$ for every $w\in W$ and $G_{\mathcal{G}%
}(\{\mathcal{G}\})=B.$ The filter $G_{\mathcal{G}}(\mathcal{G})=\mathcal{G}$
is $\mathbb{D}$-compact at $B\subset W$ by construction, so that, by
hypothesis, $R\left( \mathcal{F}\right) \times \mathcal{G}$ is $\mathbb{F}%
_{1}$-compact at $Rx\times B$ and the conclusion follows. Notice that the
inverse relation is a map $g_{\mathcal{G}}.$
\end{proof}

Theorem \ref{th:compactoidrelation} (restricted to $\mathbb{J}=\mathbb{F)}$
can be combined with Theorem \ref{th:Mperfect+range} to the effect that:

\begin{corollary}
\label{cor:prodofperfects}Let $\mathbb{D}$ and $\mathbb{M}$ be two
composable classes of filters containing principal filters. Let $%
f:X\rightarrow Y$ be a continuous surjection between two topological spaces.
The following are equivalent:

\begin{enumerate}
\item $f$ is $\mathbb{M}$-perfect with $(\mathbb{F}/\mathbb{D})$-accessible
range;

\item $f\times g$ is $\mathbb{(D\cap M)}$-perfect, for every $\mathbb{D}$%
-perfect map $g$ with $\mathbb{(D}/\mathbb{M)}$-accessible domain$;$

\item $f\times g$ is closed, for every $\mathbb{D}$-perfect map $g$ with $%
\mathbb{M}$-based domain.
\end{enumerate}
\end{corollary}

Notice that the statement corresponding to Theorem \ref%
{th:compactoidrelation} (2) is omitted in Corollary \ref{cor:prodofperfects}%
. The reason is that the hypothesis $\mathbb{M}\subset \mathbb{J}$ of
Theorem \ref{th:Mperfect+range} is in general not fullfilled so that this
statement cannot be interpreted in terms of $\mathbb{D}$-perfect maps via
Theorem \ref{th:Mperfect+range}. However, when $\mathbb{D=F},$ Theorem \ref%
{th:compactoidrelation} (2) and (4) can be interpreted properly, leading to
the following generalization of Corollary \ref{cor:Dperfectchar}.

\begin{corollary}
\label{cor:prodofperfectsD=F}Let $\mathbb{M}$ be a composable classes of
filters containing principal filters. Let $f:X\rightarrow Y$ be a continuous
surjection between two topological spaces. The following are equivalent:

\begin{enumerate}
\item $f$ is $\mathbb{M}$-perfect;

\item $f\times g$ is $\mathbb{M}$-perfect, for every perfect map $g$ with $%
\mathbb{(F}/\mathbb{M)}$-accessible domain$;$

\item $f\times g$ is closed, for every perfect map $\mathbb{M}$-based domain;

\item $f\times Id_{Y}$ is closed for every $\mathbb{M}$-based topological
space $Y.$
\end{enumerate}
\end{corollary}

The following table gathers the corresponding results. Conditions in
parenthesis are equivalent to the condition given in the same cell.

\begin{center}
\begin{tabular}{|c|c|c|c|c|}
\hline
$\mathbb{D}$ & $\mathbb{M}$ & $f\times g$ is & for every $g$ & iff $f$ is \\ 
\hline
$\mathbb{F}_{1}$ & $\mathbb{F}_{1}$ & closed & closed with & closed with \\ 
&  &  & finitely generated range & finitely generated range \\ \hline
$\mathbb{F}_{\omega }$ & $\mathbb{F}_{1}$ & closed & countably perfect with
& closed with \\ 
&  &  & finitely generated range & bisequential range \\ \hline
$\mathbb{F}_{1}$ & $\mathbb{F}_{\omega }$ & closed & closed with & countably
perfect with \\ 
&  &  & Fr\'{e}chet range & finitely generated range \\ 
&  &  & (first-countable domain) &  \\ \hline
$\mathbb{F}_{\omega }$ & $\mathbb{F}_{\omega }$ & countably perfect & 
countably perfect with & countably perfect with \\ 
&  & (closed) & strongly Fr\'{e}chet range & bisequential range \\ 
&  &  & (first-countable domain) &  \\ \hline
$\mathbb{F}_{1}$ & $\mathbb{F}$ & closed & closed & perfect with \\ 
&  &  &  & finitely generated range \\ \hline
$\mathbb{F}$ & $\mathbb{F}_{1}$ & closed & perfect with & closed \\ 
&  &  & finitely generated range &  \\ 
&  &  & (identity of finitely generated) &  \\ \hline
$\mathbb{F}_{\omega }$ & $\mathbb{F}$ & countably perfect & countably perfect
& perfect with \\ 
&  & (closed) &  & bisequential range \\ \hline
$\mathbb{F}$ & $\mathbb{F}_{\omega }$ & countably perfect & perfect with & 
countably perfect \\ 
&  & (closed) & bisequential range &  \\ 
&  &  & (identity of first-countable) &  \\ \hline
$\mathbb{F}$ & $\mathbb{F}$ & perfect & perfect & perfect \\ 
&  & (closed) & (identity map) &  \\ \hline
\end{tabular}
\end{center}

Similarily, Theorem \ref{th:compactoidrelation} (restricted to $\mathbb{J}=%
\mathbb{F)}$ can also be combined with Theorem \ref{th:Mquot+range} to the
effect that (taking again into account the restrictions applying to Theorem %
\ref{th:Mquot+range}):

\begin{corollary}
\label{cor:prodquotients} Let $\mathbb{D}$ and $\mathbb{M}$ be two
composable classes of filters containing principal filters. Let $%
f:X\rightarrow Y$ be a continuous surjection between two topological spaces.
The following are equivalent:

\begin{enumerate}
\item $f$ is $\mathbb{M}$-quotient with $(\mathbb{F}/\mathbb{D})$-accessible
range;

\item $f\times g$ is $\mathbb{(D\cap M)}$-quotient, for every $\mathbb{D}$%
-perfect map $g$ with $\mathbb{(D}/\mathbb{M)}$-accessible domain$;$

\item $f\times g$ is hereditarily quotient, for every $\mathbb{D}$-quotient
map $g$ with $\mathbb{M}$-based domain.
\end{enumerate}
\end{corollary}

\begin{corollary}
\label{cor:prodquotsD=F} Let $\mathbb{M}$ be a composable class of filters
containing principal filters. Let $f:X\rightarrow Y$ be a continuous
surjection between two topological spaces. The following are equivalent:

\begin{enumerate}
\item $f$ is $\mathbb{M}$-quotient;

\item $f\times g$ is $\mathbb{M}$-quotient, for every biquotient map $g$
with $\mathbb{(F}/\mathbb{M)}$-accessible domain$;$

\item $f\times g$ is hereditarily quotient, for every biquotient map $%
\mathbb{M}$-based domain;

\item $f\times Id_{Y}$ is hereditarily quotient for every $\mathbb{M}$-based
topological space $Y.$
\end{enumerate}
\end{corollary}

\begin{center}
\begin{tabular}{|c|c|c|c|c|}
\hline
$\mathbb{D}$ & $\mathbb{M}$ & $f\times g$ is & for every $g$ & iff $f$ is \\ 
\hline
$\mathbb{F}_{1}$ & $\mathbb{F}_{1}$ & hereditarily quotient & hereditarily
quotient with & hereditarily quotient with \\ 
&  &  & finitely generated range & finitely generated range \\ \hline
$\mathbb{F}_{\omega }$ & $\mathbb{F}_{1}$ & hereditarily quotient & 
countably biquotient with & hereditarily quotient with \\ 
&  &  & finitely generated range & bisequential range \\ \hline
$\mathbb{F}_{1}$ & $\mathbb{F}_{\omega }$ & hereditarily quotient & 
hereditarily quotient with & countably biquotient with \\ 
&  &  & Fr\'{e}chet range & finitely generated range \\ 
&  &  & (first-countable domain) &  \\ \hline
$\mathbb{F}_{\omega }$ & $\mathbb{F}_{\omega }$ & countably biquotient & 
countably biquotient with & countably biquotient with \\ 
&  & (hereditarily quotient) & strongly Fr\'{e}chet range & bisequential
range \\ 
&  &  & (first-countable domain) &  \\ \hline
$\mathbb{F}_{1}$ & $\mathbb{F}$ & hereditarily quotient & hereditarily
quotient & biquotient with \\ 
&  &  &  & finitely generated range \\ \hline
$\mathbb{F}$ & $\mathbb{F}_{1}$ & hereditarily quotient & biquotient with & 
hereditarily quotient \\ 
&  &  & finitely generated range &  \\ 
&  &  & (identity of finitely generated) &  \\ \hline
$\mathbb{F}_{\omega }$ & $\mathbb{F}$ & countably biquotient & countably
biquotient & biquotient with \\ 
&  & (hereditarily quotient) &  & bisequential range \\ \hline
$\mathbb{F}$ & $\mathbb{F}_{\omega }$ & countably biquotient & biquotient
with & countably biquotient \\ 
&  & (hereditarily quotient) & bisequential range &  \\ 
&  &  & (identity of first-countable) &  \\ \hline
$\mathbb{F}$ & $\mathbb{F}$ & biquotient & biquotient & biquotient \\ 
&  & (hereditarily quotient) & (identity map) &  \\ \hline
\end{tabular}
\end{center}

\subsection{Coreflectively modified duality}

In a series of papers \cite{DM.products}, \cite{mynard}, \cite{Mynard.survey}%
, \cite{Mynard.cap} the author developed a categorical method to deal with
topological product theorems, which relates product problems with properties
of function spaces and commutation of functors with products. Applications
of this method range from a unified treatment of a wide number of classical
results \cite{mynard}, \cite{Mynard.survey} to solutions of an old
topological problem \cite{mynard.strong} on one hand, and of a problem of
convergence theory (pertaining to Lindel\"{o}f and countably compact
convergence spaces) \cite{Mynard.comp} on the other hand. The key to
concretely apply the abstract results of \cite{mynard}, \cite{Mynard.survey}%
, \cite{Mynard.cap} is to internally characterize couples of convergences $%
(\xi ,\theta )$ (on the same underlying set) satisfying 
\begin{equation*}
\theta \times F\tau \geq G\left( \xi \times \tau \right) ,
\end{equation*}%
for every $\tau \geq H\tau $ for specific instances of concrete endofunctors 
$F,$ $G$ and $H$ of the category of convergence spaces and continuous maps.

In view of Proposition \ref{pro:local M-compactoidly}, Theorem \ref{th:main}
rephrases as follows when $A$ is a singleton.

\begin{theorem}
\label{th:commAdhD}Let $\mathbb{D}$ and $\mathbb{M}$ be two composable
classes of filters containing principal filters and let $\mathbb{J}$ be a $%
\mathbb{D}$-composable class of filters. The following are equivalent:

\begin{enumerate}
\item $\theta \geq \mathrm{Adh}_{\mathbb{J}}\mathrm{Base}_{\mathbb{D}}%
\mathrm{Adh}_{\mathbb{M}}\xi ;$

\item $\theta \times \mathrm{Base}_{\mathbb{J}}\mathrm{Adh}_{\mathbb{D}}%
\mathrm{Base}_{\mathbb{M}}S\tau \geq \mathrm{Adh}_{\mathbb{D}}\mathrm{Base%
}_{\mathbb{M}}\mathrm{Adh}_{\mathbb{M}}\left( \xi \times \tau \right) ;$

\item for every $\tau \geq \mathrm{Adh}_{\mathbb{D}}\mathrm{Base}_{\mathbb{%
M}}\tau ,$%
\begin{equation*}
\theta \times \mathrm{Base}_{\mathbb{J}}\mathrm{Adh}_{%
\mathbb{D}}\tau \geq \mathrm{Adh}_{\mathbb{D\cap M}}\left( \xi
\times \tau \right) ;
\end{equation*}

\item for every $\tau =\mathrm{Base}_{\mathbb{M}}\tau ,$%
\begin{equation*}
\theta \times \mathrm{Base}_{\mathbb{J}}\mathrm{Adh}_{%
\mathbb{D}}\tau \geq P\left( \xi \times \tau \right) .
\end{equation*}
\end{enumerate}
\end{theorem}

This generalizes \cite[Corollay 7.2 and Proposition 7.3]{mynard}
(corresponding to the case $\mathbb{J=F}$ and $\mathbb{D\subset M}$) whose
important consequences are exposed in \cite{mynard} and \cite{Mynard.survey}%
. In particular, relationships between a topological (or convergence) space
and the function spaces over it endowed with the continuous convergence (%
\footnote{%
The \emph{continuous convergence }$[\xi ,\sigma ]$ on the set of continuous
functions from $(X,\xi )$ to $(Y,\sigma )$ is the coarsest convergence
making the evaluation map jointly continuous. See \cite{BB.book} for details.%
}) can be deduced from Theorem \ref{th:main}. Beattie and Butzmann \cite%
{BB.book} call a pseudotopological space a \emph{Choquet space }and call a
space \emph{countably Choquet }if a countably based filter converges to a
point whenever all of its ultrafilter do. In other words, a convergence $\xi 
$ is countably Choquet, or in our terminology \emph{countably
pseudotopological}, if $\xi \leq \mathrm{First}S\xi .$ More generally, we
call $\mathbb{J}$\emph{-pseudotological} a convergence satisfying $\xi \leq 
\mathrm{Base}_{\mathbb{J}}S\xi $ and $\mathbb{J}$\emph{-paratopological} a
convergence satisfying $\xi \leq \mathrm{Base}_{\mathbb{J}}P_{\omega }\xi .$

Combining Theorem \ref{th:commAdhD} and \cite[Theorem 3.1]{mynard}, we get
(for $\theta =\xi $) the following new characterizations of bisequentiality,
strong and productive Fr\'{e}chetness in terms of function spaces:

\begin{corollary}
Let $\mathbb{D}$ and $\mathbb{M}$ be two composable classes of filters
containing principal filters and let $\mathbb{J}$ be a $\mathbb{D}$%
-composable class of filters. A convergence $\xi =\mathrm{Adh}_{\mathbb{M}%
}\xi $ is $\mathbb{(J}/\mathbb{D)}$-accessible if and only if $\mathrm{Base}%
_{\mathbb{J}}\mathrm{Adh}_{\mathbb{D}}\mathrm{Base}_{\mathbb{M}%
}[\xi ,\sigma ]\geq \lbrack \xi ,\sigma ]$ for every $\sigma =\mathrm{Adh}%
_{\mathbb{D}}\sigma $ (equivalently for every pretopology $\sigma $)%
$.$\newline
In particular, when $\mathbb{D=F}_{\omega }$ and $\mathbb{M=F}$:

\begin{enumerate}
\item A pseudotopology $\xi $ is bisequential if and only if the continuous
convergence $[\xi ,\sigma ]$ is a paratopology for every paratopology
(equivalently every pretopology) $\sigma ;$

\item A pseudotopology $\xi $ is productively Fr\'{e}chet if and only if $%
[\xi ,\sigma ]$ is $(\mathbb{F}_{\omega }/\mathbb{F}_{\omega })_{\#\geq }$%
-paratopological for every paratopology (equivalently every pretopology) $%
\sigma ;$

\item A pseudotopology $\xi $ is strongly Fr\'{e}chet if and only if $[\xi
,\sigma ]$ is countably paratopological for every paratopology (equivalently
every pretopology) $\sigma .$
\end{enumerate}
\end{corollary}

This is a sample example. Many others can be found in \cite{mynard} and \cite%
{Mynard.survey}.

\bibliographystyle{amsplain}

\end{document}